\documentclass[12pt,twoside,groupedaddress,titlepage,a4paper]{article}

\usepackage{amssymb}
\usepackage{amsfonts}
\usepackage{amsmath}

\def\<{\langle}
\def\>{\rangle}

\def\N{{\mathbb N}}

\def\CT{{\rm C.T.}}
\def\goth{\mathfrak}

\def\Det{\mbox{Det}}

\def\ashuff#1#2#3{
\kern 1pt \vrule height#1 \overline{\vrule height#3 width 0pt
\hskip#2} \rule{.3pt}{#1}\overline{\vrule height#3 width 0pt
\hskip#2} \rule{.3pt}{#1} \kern 1pt }

\def\Pfa{{\rm Pf\,}}

\def\S{\goth S}

\title{Hyperdeterminantal calculations of   Selberg's and
Aomoto's integrals}

\date{\small Dedicated to the Memory of Brian G. Wybourne}

\author{Jean-Gabriel Luque
\thanks{IGM, Universit\'e de Marne-la-Vall\'ee, 77454 Marne-la-Vall\'ee Cedex 2, France ({\tt luque@univ-mlv.fr})}
\and Jean-Yves Thibon
\thanks{IGM, Universit\'e de Marne-la-Vall\'ee, 77454 Marne-la-Vall\'ee Cedex 2, France ({\tt jyt@univ-mlv.fr})}}

\begin{document}

\maketitle

\begin{abstract}
The hyperdeteminants considered here
are the simplest analogues of determinants for higher rank tensors
which have been defined by Cayley, and apply only to
tensors with an even number of indices. We have shown in a previous
article that the calculation of certain multidimensional integrals could
be reduced to the evaluation of hyperdeterminants of Hankel type.
Here, we carry out this computation by purely algebraic means in
the cases of Selberg's and Aomoto's integrals.
\end{abstract}

\section{Introduction}
Selberg's integral and its generalizations are multiple integrals
encountered in many fields such as  mathematical statistics,
random matrices,
statistical mechanics, special functions, integrable systems,
and can be interpreted as generalized moments,
correlation functions or partition functions.

Such integrals appear naturally in the computation of
the normalization constants for various random matrix
ensembles. For example the
normalization of the eigenvalue probability density function for
Gaussian, circular, Laguerre, Jacobi and Cauchy ensembles are  limit cases
of Selberg integrals (see \cite{Me,Fo1,Fo2}).
Generalizations appear in the calculation of the partition functions
of log-potential Coulomb systems  \cite{Fo3}.
Such multiple moment integrals
are also connected with the problem of expanding
even powers of the Vandermonde determinant
\begin{equation}
\Delta(x)=
\det(x_i^{j-1}) =\prod_{1\le i<ji \le n}(x_j-x_i)
\end{equation}
in various bases
of symmetric functions.
Expansion of these powers in terms of  Schur functions is still an
open problem (see \cite{STW,Wy1,KTW}), whose solution
would find important applications  to
calculations related to some aspects of the fractional quantum Hall effect (see
\cite{STW,Wy2}) described by Laughlin's wave function \cite{Lau}
\begin{equation}
\Psi^m(z_1,\dots,z_n)=\prod_{i<j}^N(z_i-z_j)^{2m+1}\exp\left\{-\frac12\sum_{i=1}^N
|z_i|^2\right\}
\end{equation}
The same expansion would allow one to compute many more interesting
 integrals (see
\cite{Fo1,Fo3,LT1,LT2}).

Selberg's integral is a generalization of Euler's
Beta integral, introduced  by himself in 1944 (see
\cite{Sel}), in view of proving   a generalization of a theorem
of Gelfond.
Selberg  found that his integral could be
written as a product of quotients of gamma functions
\begin{eqnarray}
{\cal
S}_{n}(a,b;\gamma)&=&\int_0^1\cdots\int_0^1|\Delta(x)|^{2\gamma}
\prod_{i=1}^nx_i^{a-1}(1-x_i)^{b-1}dx_i\nonumber\\
&=&\prod_{j=0}^{n-1}
{\Gamma(a+j\gamma)\Gamma(b+j\gamma)\Gamma((j+1)\gamma+1)
\over
\Gamma(a+b+(n+j-1)\gamma)\Gamma(\gamma+1)}.
\label{Selberg}\end{eqnarray}
This identity was subsequently proved in many
different ways by other authors, such as
Anderson \cite{An}, Aomoto \cite{Ao1} and Dotsenko and Fateev
\cite{DF}. These proofs are all interesting, as they involve
generalizations of  Selberg's integral or techniques which can be
useful in other contexts. However,  Selberg's original proof
has the specificity that it is almost elementary, in the sense that he
 computed the integral in the case where $\gamma$ is an integer $k$,
using only arguments on the degrees of the polynomials appearing after
integration of the terms obtained by expansion of $\Delta(x)^{2k}$,
and a very simple change of
variables. Then, he concluded by involving a classical theorem
of Carlson
to extend it to complex values of $\gamma$ (see \cite{Me} or \cite{Ti}).
As we shall see below, the  elementary part of this proof comes out naturally if we
consider this object for $\gamma=k\in\N$ not as an integral but as
a special case of a
 hyperdeterminant, a notion originally due to Cayley
but whose only available textbook presentations seem to be Refs.
\cite{Sok1, Sok2} (in Russian).

The hyperdeterminant is a  polynomial invariant
of tensors of even order ({\it i.e.}, with an even number of indices).
In the case of matrices, its definition gives back
the classical determinant.

We have explained in \cite{LT2} that Selberg's integral could be
rewritten as a hyperdeterminant of Hankel type ({\it i.e.},  whose
entries depend only on the sum of the indices) built from
the moments of the beta distribution. In this paper, we give a complete
translation of
Selberg's proof in the hyperdeterminantal formalism, using only very
simple tools. The hyperdeterminantal aspects of  multiple
integral computations look very promising. At present, we
know only a few tools to handle such  polynomials, but this is
already sufficient to deal with some interesting cases. To illustrate this point,
we sketch a hyperdeterminantal proof of Aomoto's integral
\begin{eqnarray}\label{Aomoto}
{\cal
A}_n^{a,b;k}(y)&=&\int_0^1\cdots\int_0^1\Delta(x)^{2k}\prod_{i=1}^n(y-x_i)x_i^{a-1}(1-x_i)^{b-1}dx_1\cdots dx_n\\
&=&(-2)^{-n}S_{n}(a,b;k)P_n^{{a\over k}-1,{b\over k}-1}(1-2y)
\end{eqnarray}
where $P_n^{a,b}$ are the monic Jacobi polynomials.
Aomoto has proved this identity  in \cite{Ao1} by means of
differential equations.
The main interest of   Aomoto's original proof is its relation with
the Calogero model. Barsky and Carpentier have given in \cite{BC} a
different proof, based on Anderson's method.

\section{Elementary hyperdeterminant calculus}

\subsection{Definition}

We shall find it convenient to use a set of $n$
Grassmann (that is, anticommutative,
or fermionic) variables
$\eta=\{\eta_1,\dots,\eta_i,\dots\}$.
A tensor of order $k$ and
dimension $n$ , that is, an element of $V^{\otimes k}$, where $V$
is an $n$ dimensional vector space, will be represented as
\begin{equation}
{\bf M}=\sum_{1\leq
i_1,\cdots,i_{k}\leq n}M_{i_1\cdots i_{k}}\eta_{i_1}\otimes\cdots\otimes\eta_{i_k}.
\end{equation}
That is, we identify $V$ with the vector space spanned by $\eta_1$,
$\dots,$ $\eta_n$.
The array $(M_{i_1 \dots i_k})$ will be called a hypermatrix.
The hyperdeterminant $\Det({\bf M})$ of $\bf M$ can be compactly defined as the  coefficient of
$(\eta_{1}\cdots\eta_n)^{\otimes k}$ in ${\bf M}^n\over n!$:
\begin{equation}\label{defhypdet}
{\bf M}^n=n!\,\Det({\bf M})(\eta_{1}\cdots\eta_n)^{\otimes k}
\end{equation}
where the power ${\bf M}^n$ is evaluated in the tensor product of
$k$ copies of the Grassmann algebra over the  $n$ distinct
$\eta_i$'s, so that the only surviving terms are those involving
$(\eta_{1}\cdots\eta_n)^{\otimes k}$.

For convenience, we will use the notation
\begin{equation}
\Det_{k}(M_{i_1 \cdots i_{k}})_1^{n}=\Det\left({\bf M}\right).
\end{equation}
Note first that if $\bf M$ is a tensor of odd order $k$, one has
$\Det({\bf M})=0$. If $\bf M$ is a tensor of even order $2k$ and dimension $n$, one
has
\begin{equation}\label{sumdet}
\Det({\bf M})={1\over
n!}\sum_{\sigma_1,\cdots,\sigma_{2k}\in\S_n}\epsilon(\sigma_1)\cdots\epsilon(\sigma_{2k})\prod_{i=1}^nM_{\sigma_1(i) \dots
\sigma_{2k}(i)}.
\end{equation}
Setting $k=1$, we recover the classical definition of the determinant
for matrices.

\subsection{The invariance property}
Definition (\ref{defhypdet}) is very useful for proving properties of hyperdeterminants.
 Let us consider the
natural action of the group $GL(V)^{\times 2k}$ on $V^{\otimes 2k}$
\begin{equation}
(g^{(1)}\otimes\cdots\otimes g^{(2k)})\cdot{\bf v}_1\otimes\cdots\otimes{\bf
v}_{2k}=(g ^{(1)}{\bf v}_1)\otimes\cdots\otimes(g ^{(2k)}{\bf v}_{2k})
\end{equation}
Then, for a general ${\bf M}\in V^{\otimes 2k}$,
\begin{eqnarray}\label{eq10}
(g^{(1)}\otimes\cdots\otimes g^{(2k)})\cdot{\bf M}=\displaystyle\sum_{1\leq j_1,\cdots,j_{2k}\leq n}
M_{j_1 \cdots j_{2k}}\times\nonumber\\
\displaystyle\times \sum_{i_1=1}^ng^{(1)}_{i_1 j_1}\eta_{j_1}\otimes\cdots\otimes\sum_{i_{2k}=1}^ng^{(2k)}_{i_{2k} j_{2k}}
\eta_{j_{2k}}.
\end{eqnarray}
Hence,
\begin{equation}\label{invariance}
\Det\left((g^{(1)}\otimes\cdots\otimes g^{(2k)})\cdot{\bf
M}\right)=\det\left(g^{(1)}\right)\cdots\det\left(g^{(2k)}\right)\Det({\bf
M}).
\end{equation}

\subsection{Minor summation formula}
Consider two tensors
\begin{equation}
{\bf M}=\sum_{1,\leq i_1 \ldots, i_{2k}\leq
n}M_{i_1 \dots  i_{2k}}\eta_{i_1}\otimes\cdots\otimes\eta_{i_{2k}}
\end{equation}
and
\begin{equation}
{\bf N}=\sum_{1\leq i_1,\dots,i_{2k}\leq
n}N_{i_1 \dots i_{2k}}\eta_{i_1}\otimes\cdots\otimes\eta_{i_{2k}}.
\end{equation}
Remarking that ${\bf MN}-{\bf NM}=0$
(in the $2k$th tensor power of the Grassmann algebra), one has
\begin{equation}\label{binome}
({\bf M}+{\bf N})^n=\sum_{i=0}^n\left(n\atop i\right){\bf M}^i{\bf
N}^{n-i}.
\end{equation}
If $I=(I_1,\cdots,I_{2k})$ is a $2k$-tuple of subsets of
$\{1,\dots,n\}$, we will denote by ${\bf M}[I]$ the tensor
\begin{equation}
{\bf M}[I]=\sum_{i_1\in I_1,\dots,i_{2k}\in
I_{2k}}M_{i_1 \cdots i_{2k}}\eta_{i_1}\otimes\cdots\otimes\eta_{i_{2k}}
\end{equation}
and by $\eta_I$ the product
\begin{equation}
\eta_I=\overrightarrow{\prod_{i_1\in
I_1}}\eta_{i_1}\otimes\cdots\otimes\overrightarrow{\prod_{i_{2k}\in
I_{2k}}}\eta_{i_{2k}}
\end{equation}
where the symbol $\overrightarrow{\prod_i}$ denotes the
 product  taken in increasing order of the subscripts.

Let ${\goth C}_{n,k}^r$ be the set of  pairs $(I,J)$ of
$2k$-tuples $I=(I_1,\cdots,I_{2k})$ and $J=(J_1,\cdots,J_{2k})$
such that for each $s\in\{1,\dots,2k\}$, $(I_s,J_s)$ is a
partition of $\{1,\dots,n\}$ into two blocks of sizes $r$ and $n-r$.
For each $(I,J)\in{\goth C}_{n,k}^r$ we will denote by
$\epsilon(I,J)$ the product of the signs of the permutations
$\sigma_s=(i_1,\dots,i_r,j_1,\dots,j_n)$ (written as a word) where
$I_s=\{i_1<\cdots<i_r\}$ and $J_s=\{j_1<\cdots<j_{n-r}\}$.

We can now write, for each $p\leq n$,
\begin{equation}
{\bf M}^p=p!\sum_{I}\Det({\bf M}[I])\eta_I
\end{equation}
where the sum is over the $2k$-tuples $I=(I_1,\cdots,I_{2k})$ of
subsets of $\{1,\dots,n\}$ of cardinality $p$. Hence, by
(\ref{binome}),
\begin{equation}\label{sumtensor}
\Det({\bf M}+{\bf N})=\sum_{r=0}^n\sum_{(I,J)\in{\goth C}^r_{n,k}}\epsilon(I,J)\Det({\bf
M}[I])\Det({\bf N}[J])
\end{equation}
Now, if we set
\begin{equation}
{\bf M}'=\sum_{1\leq i_2,\dots,i_{2k}\leq
n}M_{1,i_2,\dots,i_{2k}}\eta_1\otimes\eta_{i_2}\otimes\cdots\otimes
\eta_{i_{2k}}
\end{equation}
and
\begin{equation}
{\bf M}''=\sum_{i_1=2}^n\sum_{1\leq i_2,\dots,i_{2k}\leq
n}M_{i_1,i_2,\dots,i_{2k}}\eta_{i_1}\otimes\cdots\otimes
\eta_{i_{2k}},
\end{equation}
remarking that ${\bf M}={\bf M}'+{\bf M}''$ we obtain from
(\ref{sumtensor}),
\begin{eqnarray}\label{dev1ind}
\Det({\bf M})=&\displaystyle\sum_{I=(i_1=1\leq i_2,i_3,\dots,i_{2k}\leq
n)}{\rm sign}(I)M_{1,i_2,\dots,i_{2k}}\Det({\bf M}[\overline I])
\end{eqnarray}
where ${\rm sign}(I)=(-1)^{i_1+i_2+\cdots+i_{2k}}$ and $\overline
I=(\{1,\dots,n\}-i_1,\dots,\{1,\dots,n\}-i_{2k})$.\\
Equations (\ref{sumtensor}) and (\ref{dev1ind}) can be found in
\cite{Bar,Sok1,Sok2}.

\section{Hankel hypermatrices}
For a sequence $I=(i_1,\dots,i_{2k})$, we will set
$|I|=i_1+\cdots+i_{2k}$ and
$\eta_I=\eta_{i_1}\otimes\cdots\otimes\eta_{i_{2k}}$.

\subsection{Heine's integrals for Hankel determinants}

Let $\mu$ be a measure
on the real line,  and $c_n=\int x^nd\mu(x)$ be its moments.
Heine's integrals below are evaluated as
Hankel determinants whose entries are the moments $c_n$ (see
\cite{He})
\begin{equation}
\int\cdots\int\Delta(x)^2d\mu(x_1)\cdots d\mu(x_n)=n!\det(c_{i+j})_0^{n-1}.
\end{equation}
It is easy to give a direct
proof of this equality, but
it can also be seen as a very special case of  de Bruijn's
integral
\begin{equation}\label{fdB2}
\mathop{\int\cdots\int}_{a\leq x_1< \cdots< x_{n}\leq b} \det
\left(\phi_i(x_j)|\psi_i(x_j)\right)dx_1\cdots
dx_n=\Pfa\left(Q_{ij}\right)_{1\leq i,j\leq 2n}
\end{equation}
where $Q_{ij}=\int_a^b[\phi_i(x)\psi_j(x)-\phi_j(x)\psi_i(x)]dx$
 and $\left(\phi_i(x_j)|\psi_i(x_j)\right)$ denotes the matrix
whose $i$th row is
$[\phi_i(x_1),\psi_i(x_1),\phi_i(x_2),\psi_i(x_2),\ldots,\phi_i(x_n),\psi_i(x_n)]$
and $\Pfa(M)$ denotes the Pfaffian of the matrix $M$. In
\cite{LT1}, we have obtained a generalization of this identity, involving
$2kn$ functions on the left hand side and a hyperpfaffian (see
\cite{Bar}) on the right hand side. In the case where the
determinant is an even power of the Vandermonde determinant, we
obtain a generalization of Heine's theorem turning multiple integrals into
Hankel hyperdeterminants. This generalization is discussed in
\cite{LT2} and will be recalled in the following subsection.

\subsection{Generalization: Hankel hyperdeterminants}

We consider here the integral
\begin{equation}
{\cal I}_{n,k}^\mu=\int\cdots\int\Delta(x_1,\cdots,x_n)^{2k}d\mu(x_1)\cdots
d\mu(x_n)
\end{equation}
Expanding the even power of the Vandermonde determinant, this
integral can be expressed as a Hankel hyperdeterminant whose entries
are the moments $c_n=\int x^nd\mu(x)$:
\begin{equation}\label{Heine}
{\cal I}_{n,k}^\mu=n!\Det\left(c_{|I|}\right)_0^{n-1}.
\end{equation}
Consider now the  following two tensors
\begin{equation}
{\bf S}_{n}(a,b;k)=\sum_{0\leq i_1,i_2,\dots,i_{2k}\leq n-1}{\rm
B}(a+|I|,b)\eta_{I}
\end{equation}
and
\begin{equation}
{\bf A}_{n}^{a,b;k}(y)=\sum_{0\leq i_1,i_2,\dots,i_{2k}\leq
n-1}{\rm
B}(a+|I|,b)\cdot \left(y-{a+|I|\over a+b+|I|}\right)\cdot \eta_I,
\end{equation}
where ${\rm B}(a,b)$ denotes the Beta function ${\rm
B}(a,b)=\int_0^1x^{a-1}(1-x)^{b-1}dx$.
Remarking that the ${\rm B}(a+n,b)$ are the moments of the Beta
distribution, one obtains after applying (\ref{Heine})
\begin{equation}\label{SelbHyp}
{\cal S}_n(a,b;k)=n!\Det\left({\bf S}_n(a,b;k)\right)
\end{equation}
and in the same way
\begin{equation}\label{AomotoHyp}
{\cal A}_n^{a,b;k}(y)=n!\Det\left({\bf A}_n^{a,b;k}(y)\right).
\end{equation}

\subsection{Expansion of Hankel hyperdeterminants}

In this section we consider the hyperdeterminant of a general Hankel
tensor
\begin{equation}
{\cal H}_{n,k}=\Det_{2k}\left(X_{i_1+\cdots+i_{2k}}\right)_0^{n-1}.
\end{equation}
Expanding this as a polynomial in the $X_i$, one finds
\begin{equation}\label{devhank}
{\cal
H}_{n,k}=\sum_{\lambda}c_{\lambda}^{n,k}X_{\lambda_1}\cdots
X_{\lambda_n}
\end{equation}
where the sum is over the $n$-tuples
$\lambda=(\lambda_1,\cdots,\lambda_n)$ such that
$|\lambda|=kn(n-1)$ and such that for each $i\in\{1,\cdots,n\}$, $(i-1)k\leq
\lambda_i\leq k(n+i-2)$.

This expansion is derived in \cite{LT2}.
A special case appears in  Selberg's  proof (see
\cite{Me,Sel}).

Let  ${\goth d}_{n,k}$ be the
coefficient of $X_{k(n-1)}^n$ in (\ref{devhank}):
\begin{equation}\label{nabla}
{\goth d}_{n,k}=c_{(k(n-1))^n}.
\end{equation}
If we set $X_i=0$ for $i>k(n-1)$, we obtain from (\ref{devhank}) that
\begin{equation}\label{detnabla}
\Det_{2k}\left(X_{i_1+\cdots+i_{2k}}\right)_0^{n-1}={\goth d}_{n,k}X_{k(n-1)}^n.
\end{equation}
A closed form for  ${\goth d}_{n,k}$ follows from
a well-known constant term identity  (the
``Dyson conjecture'')\footnote{This was conjectured by  Dyson \cite{Dy} in 1962.
I. J. Good gave, in 1970, an elegant elementary  proof involving
only Lagrange interpolation \cite{Go}.}
\begin{equation}\label{Dyson}
{\cal C}_{n,k}=\CT\left\{\prod_{i,j=1\atop i\neq j }^n
\left(1-{x_i\over
x_j}\right)^{a_i}\right\}=\left(a_1+\cdots+a_n\atop
a_1,\cdots,a_n\right),
\end{equation}
where C.T. means ``constant term''.
On the other hand,
\begin{equation}
\prod_{i,j=1\atop i\neq
j}^n\left(1-{x_i\over
x_j}\right)^k=(-1)^k\Delta(x)^{2k}\prod_{i=1}^n\frac1{x_i^{k(n-1)}}.
\end{equation}
Expanding the power of the Vandermonde determinant, one finds that
for $a_i\equiv k$,
\begin{eqnarray}
{\cal C}_{n,k}&=&(-1)^k\CT\left\{\sum_{\sigma_1,\cdots,\sigma_{2k}\in\S_n}
\epsilon(\sigma_1)\cdots\epsilon(\sigma_{2k})\prod_{i=1}^n{x_i^{\sigma_1(i)+\cdots+\sigma_{2k}(i)-2k-k(n-1)}}\right\}\nonumber\\
&=&(-1)^kn!\Det_{2k}\left(\delta_{i_1+\cdots+i_{2k},n}\right)_0^{n-1}
\end{eqnarray}
where $\delta_{i,j}=1$ if $i=j$ and $0$ otherwise
(Kronecker symbol).
Using (\ref{detnabla}), we obtain
\begin{equation}\label{ct}
{\cal C}_{n,k}=(-1)^kn!{\goth d}_{n,k}.
\end{equation}
Hence,
\begin{equation}\label{ctnabla}
{\goth d}_{n,k}=(-1)^k{1\over n!}\left(kn\atop
k,\cdots,k\right).
\end{equation}

\subsection{Minors of a Hankel hyperdeterminant}

We consider here a Hankel hyperdeterminant
\begin{equation}
{\bf
M}=\sum_{I}M_{|I|}\eta_{i_1}\otimes\cdots\otimes
\eta_{i_{2k}}
\end{equation}
and a family of  $2k$  subsets $J=(J^1,\dots,J^{2k})$ of a fixed
size $m\leq n$  with $J^l=\{j^l_0\leq j^l_1\leq\cdots\leq
j^l_{m}\}\subset \{0,\dots,n-1\}$. We want to expand the
polynomial $M[J]^m$, where
\begin{equation}
{\bf M}[J]=\sum_{I}M_{j^1_{i_1}+\cdots+j_{i_{2k}}^{2k}}\eta_{i_1}\otimes\cdots\otimes
\eta_{i_{2k}}.
\end{equation}
One has
\begin{equation}
\Det({\bf M}[J])=\frac1{n!}\sum_{\sigma_1,\dots,\sigma_{2k}\in\S_m}
\epsilon(\sigma_1)\cdots \epsilon(\sigma_{2k})
\prod_{l=1}^m
M_{j^1_{\sigma_1(l)}+\cdots+j_{\sigma_{2k}(l)}^{2k}}.
\end{equation}
But for each $2k$-tuple $(\sigma_1,\dots,\sigma_{2k})$, one has
\begin{eqnarray}\prod_{l=1}^m
M_{j^1_{\sigma_1(l)}+\cdots+j_{\sigma_{2k}(l)}^{2k}}=\prod_{l=1}^mM_{\sigma_1(l)+\cdots+
\sigma_{2k}(l)-2k+r^1_{\sigma_1(l)}+\cdots+r^{2k}_{\sigma_{2k}(l)}}
\end{eqnarray}
with $0\leq r^i_l\leq n-m$. Hence, by a  computation
similar to that of a complete Hankel hyperdeterminant,
\begin{eqnarray}\prod_{l=1}^m
M_{j^1_{\sigma_1(l)}+\cdots+j_{\sigma_{2k}(l)}^{2k}}=\prod_{l=1}^mM_{\lambda_i+
s_i}
\end{eqnarray}
with $k(i-1)\leq \lambda_i\leq k(m+i-2)$ and $0\leq s_i\leq
2k(n-m)$. And
\begin{equation}\label{expmin}
\Det(M[J])
=
\sum_{\tilde\lambda}
c^{n,k;J}_{\tilde\lambda}
\prod_{i=1}^m\tilde
M_{\tilde\lambda_i} \end{equation} where
$k(i-1)\leq\tilde\lambda_i\leq k(2n-m+i-2)$.

\section{Selberg's integral}

\subsection{From the Hankel form to the symmetric form}

Selberg's integral admits another hyperdeterminantal
representation which can be obtained directly
from the previous one  (without manipulating the
integral). It suffices to remark the following property
\begin{equation}\label{Beta1}
\sum_{i=0}^n(-1)^i\left(n\atop i\right){\rm B}(a+i,b)={\rm
B}(a,b+n).\end{equation}
Hence, we have from (ref{eq10})
\begin{equation}
{\bf S}_n(a,b;k)=\displaystyle(\overbrace{{\rm I},\cdots,{\rm I}}^{
k},\overbrace{{ g},\cdots,{ g}}^{ k})\cdot {\bf S}^{\it Sym}_n(a,b;k)
\end{equation}
where
${\rm I}$ is the identity matrix,  ${g}$ is the $n\times n$ matrix
\begin{equation}
{g}=\sum_{1\leq i,j\leq n}(-1)^{j}\left(i\atop
j\right)\eta_i\otimes\eta_j,
\end{equation}
and
\begin{equation}
{\bf S}^{\it
Sym}_n(a,b;k)=\sum_{I,J\in\{0,\dots,n-1\}^{k}}
{\rm B}(a+|I|,b+|J|)\eta_{I}\otimes\eta_{J}
\end{equation}
It follows from (\ref{invariance}) that
\begin{equation}
\det({g})^k\Det\left({\bf S}_n^{Sym}(a,b;k)\right)=\Det\left({\bf
S}_n(a,b;k)\right)
\end{equation}
and, since
\begin{equation}
\det({g})=(-1)^{n(n-1)\over 2},
\end{equation}
we find
\begin{equation}\label{symmetry}
\Det\left({\bf
S}_n^{Sym}(a,b;k)\right)=(-1)^{kn(n-1)\over 2}\Det\left({\bf
S}_n(a,b;k)\right)
\end{equation}
without the help of the integral representation.
\subsection{Selberg's original proof}
Each step of  Selberg's proof can be viewed as a
manipulation on a hyperdeterminant (see \cite{Me} for example). It can be
divided into two parts. In the first part, Selberg proved that
\begin{equation}\label{eqs1}
{\cal S}_n(a,b;k)=c_{n,k}\prod_{i=0}^{n-1}{\Gamma(a+jk)\Gamma(b+jk)\over
\Gamma(a+b+(n+j-1)k)}.
\end{equation}
where $c_{n,k}$ does not depend on $a$ and $b$. In this proof,
Selberg started from the expansion of the Vandermonde determinant.
In term of hyperdeterminants, it is equivalent to the expansion
(\ref{SelbHyp}).
For our purpose,  equation (\ref{eqs1}) reads
\begin{equation}\label{eqd1}
\Det({\bf S}_n(a,b;k))=\alpha_{n,k}\prod_{i=0}^{n-1}{\Gamma(a+jk)\Gamma(b+jk)\over
\Gamma(a+b+(n+j-1)k)}
\end{equation}
where $\alpha_{n,k}=n!c_{n,k}$. To obtain it, we start from
formula (\ref{devhank})
\begin{equation}
\Det({\bf S}_n(a,b;k))=\sum_{\lambda}c_{\lambda}^{n,k}{\rm B}(a+{\lambda_1},b)\cdots
{\rm B}(a+ {\lambda_n},b)
\end{equation}
where the sum is over all  $n$-tuples
$\lambda=(\lambda_1,\cdots,\lambda_n)$ such that
 $|\lambda|=kn(n-1)$ and for each $i\in\{1,\cdots,n\}$, $(i-1)k\leq
\lambda_i\leq k(n+i-2)$. Then, extracting common factors,  we
arrive at an expression of the form
\begin{equation}
\Det({\bf S}_n(a,b;k))={Q(a,b)\over R(b)}
\prod_{i=0}^{n-1}{\Gamma(a+jk)\Gamma(b+jk)\over
\Gamma(a+b+(n+j-1)k)}
\end{equation}
where $Q(a,b)$ is a polynomial of degree at most $kn(n-1)/2$
and $R(b)$ is a polynomial of degree $kn(n-1)/2$. But
(\ref{symmetry}) implies that
\begin{equation}
\Det({\bf S}_n(a,b;k))=\Det({\bf S}_n(b,a;k))
\end{equation}
and then
\begin{equation}
{Q(a,b)\over R(b)}={Q(b,a)\over R(a)}
\end{equation}
which implies that $\alpha_{n,k}={Q(a,b)\over R(b)}$ is
independent of $a$ and $b$.\\

The second part of  Selberg's proof consists in establishing  a
recurrence relation for  $\alpha_{n,k}$:
\begin{equation}\label{recc}
\alpha_{n,k}={(nk)!\over n!}\alpha_{n-1,k}.
\end{equation}
This can be found from the limit
\begin{equation}
L_{b;k,n}=\lim_{a\rightarrow 0^+}a\,\Det({\bf S}_n(a,b;k))
\end{equation}
Hence, one has
\begin{equation}
L_{b;k,n}=\lim_{a\rightarrow 0^+}\Det({\bf T}_a)
\end{equation}
where ${\bf T}_a$ is the tensor defined by
\begin{eqnarray}
{\bf T}_a=\sum_{I\in\{0\}\times\{0,\dots,n-1\}^{2k-1}}a{\rm
B}(a+|I|,b)\eta_I+\nonumber\\+
\sum_{I\in\{1,\cdots,n-1\}\times\{0,\dots,n-1\}^{2k-1}}{\rm
B}(a+I,b)\eta_I
\end{eqnarray}
Remarking that
\begin{equation}
\lim_{a\rightarrow0^+}a\,{\rm B}(a+n,b)=\delta_{n,0}
\end{equation}
where $\delta_{i,j}$ is the Kronecker symbol, one obtains
\begin{equation}
L_{b;k,n}=\Det({\bf T})
\end{equation}
where $\bf T$ is the tensor
\begin{eqnarray}
{\bf
T}=\eta_0\otimes\cdots\otimes\eta_0+\sum_{I\in\{1,\cdots,n-1\}\times\{0,\dots,n-1\}^{2k-1}}{\rm
B}(|I|,b)\eta_{I}
\end{eqnarray}
Hence, expanding by (\ref{dev1ind}) with respect to the first subscript, one has
\begin{equation}
L_{b;k,n}=\Det({\bf S}_{n-1}(1,b;k)).
\end{equation}
Together with (\ref{eqd1}) this equality gives the recurrence
relation (\ref{recc}) and proves  Selberg's identity.

\subsection{An alternative ending for Selberg's proof}

We shall now give a simpler proof of Selberg's identity by reducing it to the
Dyson conjecture, which, as we have already seen, is by now a
familiar and  elementary
statement involving no more than  Lagrange interpolation.
Using the functional equation of the $\Gamma$ function and the
linearity properties of the hyperdeterminant, we
can recast
the Hankel hyperdeterminantal expression of Selberg's
integral  in the form
\begin{equation}
\Det\left({\bf S}_n(a,b;k)\right)={\rm
B}(a,b)^n\Det_{2k}\left({(a)_{|I|}\over(a+b)_{|I|}}\right)_0^{n-1}.
\end{equation}
where $(a)_n=a(a+1)\cdots (a+n-1)$ denotes the Pochhammer symbol.
The result will follow if we can
obtain a closed form for the hyperdeterminant
\begin{equation}
{\cal
D}_{n,k}(a,b)=\Det_{2k}\left({(a)_{|I|}\over(b)_{|I|}}\right).
\end{equation}
We start as in the previous proof, obtaining an analog of (\ref{eqd1})
\begin{equation}
{\cal
D}_{k,n}(a,b)=\alpha'_{n,k}\prod_{m=1}^n{(a)_{k(m-1)}(b-a)_{k(m-1)}\over(b)_{k(n+m-1)}}\end{equation}
where the constant $\alpha'_{n,k}$ is independent of $a$ and
$b$.
Now, for every $a'$
\begin{equation}
{\cal D}_{k,n}(a,b)=\prod_{m=1}^n{(a)_{k(m-1)}(b-a)_{k(m-1)}\over(b)_{k(n+m-1)}(a')_{k(m-1)}
}\Det_{2k}\left((a')_{|I|}\right)_0^{n-1}.
\end{equation}
This means that Selberg's integral can be deduced from the
knowledge of $\Det_{2k}((a')_{|I|})_0^{n-1}$ for any particular
choice of
$a'$. If we set $a'=-k(n-1)$, one has by
(\ref{detnabla})
\begin{equation}
\Det_{2k}\left((-k(n-1))_{|I|}\right)_0^{n-1}=(-1)^{kn(n-1)}{\goth d}_{n,k}(k(n-1))!^n.
\end{equation}
But ${\goth d}_{n,k}$ has already been evaluated by means of the Dyson
conjecture (\ref{ctnabla}) and the closed form of  Selberg's
integral follows.

\section{Aomoto's integral}

\subsection{Another hyperdeterminantal representation}
On the integral representation (\ref{Aomoto}), we obtain easily the
following identity through the substitution $x_i\rightarrow 1-x_i$
\begin{equation}
{\cal A}_n^{a,b;k}(y)=(-1)^n{\cal A}_n^{b,a;k}(1-y).
\end{equation}
There is also a simple hyperdeterminantal proof of this assertion. We
consider the tensor
\begin{eqnarray}
{\bf A}_n^{b,a;k}(1-y)=\displaystyle\sum_{I\in\{0,\dots,n-1\}^{2k}}\left(1-y-{\displaystyle b+|I|\over
\displaystyle a+b+|I|}\right)
{\rm B}(b+|I|,b)\eta_{I}
\end{eqnarray}
and the action of $2k$ copies of the matrix ${g}$ defined before.
Using equality (\ref{Beta1}), one has
\begin{equation}
({g},\dots,{g}){\bf A}_n^{b,a;k}(1-y)=-{\bf A}_n^{a,b;k}(y)
\end{equation}
Hence,
\begin{equation}\label{chvar}
\Det({\bf A}_n^{a,b;k}(y))=(-1)^n\Det({\bf A}_n^{b,a;k}(1-y)).
\end{equation}

\subsection{Minors of Selberg's hyperdeterminant}

To prove Aomoto's identity, we need some preliminary results on Selberg's
hyperdeterminant.
Let us consider a family of $2k$ subsets
$I=(I^1,\dots,I^{2k})$ with $I^{j}=\{i_0^j\leq i_1^j,\dots,\leq
i_{m-1}^j \}\subset\{1,\cdots,n\}$ and the associated sub-tensor
of  Selberg's Hankel tensor:
\begin{equation}
{\bf S}_n(a,b;k)[I]=\sum_{0\leq j_1,\dots,j_{2k}\leq m-1}{\rm B}
(a+i_{j_1}^1+\cdots+i_{j_{2k}}^{2k},b)
\eta_{i_{j_1}^1}\otimes\cdots\otimes\eta_{i_{j_k}^{2k}}.
\end{equation}
One can write by (\ref{expmin})
\begin{equation}
\Det({\bf S}_n(a,b;k)[I])=\sum_{\lambda}
c_\lambda^{n,k;I}\prod_{i=1}^m{\rm B}(a+\lambda_i)
\end{equation}
where $k(i-1)\leq\lambda_i\leq k(2n-m+i-2)$. Hence,
\begin{eqnarray}\label{minor1}
\Det({\bf
S}_n(a,b;k)[I])=&\displaystyle\prod_{i=1}^m{\Gamma(a+k(i-1))\Gamma(b)\over\Gamma(a+b+k(2n-m+i-2))}{{\goth
p}_I(a,b)}
\nonumber\\&
=\displaystyle\prod_{i=1}^m{\Gamma(a+k(i-1))\Gamma(b)\over\Gamma(a+b+k(2n-i-1))}{{\goth
p}_I(a,b)}
\end{eqnarray}
where  ${\goth p}_I(a,b)$ is a polynomial. In the same way, we
will use the following identity
\begin{eqnarray}\label{minor2}
\Det({\bf S}_n(a+1,b;k)[I])=\prod_{i=1}^n{1\over
a+b+k(2n-i-1)}\prod_{i=1}^m(a+k(i-1))\times\nonumber
\\
\times\prod_{i=m+1}^n(a+b+k(2n-i-1))
\prod_{i=1}^m{\Gamma(a+k(i-1))\Gamma(b)\over\Gamma(a+b+k(2n-i-1))}{{\goth
p}_I(a+1,b)}.
\end{eqnarray}
which follows immediately from (\ref{minor1}).

\subsection{A proof of Aomoto's identity}

We start from the hyperdeterminantal representation and we remark that
\begin{equation}
{\bf A}^{a,b;k}_n(y)={\bf S}(a,b;k)y-{\bf S}(a+1,b;k).
\end{equation}
From (\ref{sumtensor}), we obtain
\begin{eqnarray}
\Det({\bf
A}_n^{a,b;k}(y))=\sum_{r=0}^n(-1)^{n-r}y^r\times\nonumber\\\times\sum_{(I,J)\in{\goth
C}^r_{n,k}}\epsilon(I,J)\Det({\bf S}(a,b;k)[I])\Det({\bf
S}(a+1,b;k)[J]).
\end{eqnarray}
From (\ref{minor1}) and (\ref{minor2}), one has
\begin{eqnarray}\label{Ao1}
\Det({\bf A}_n^{a,b;k}(y))=\prod_{i=1}^n{1\over
a+b+k(2n-i-1)}\times\nonumber\\
\times\sum_{r=0}^ny^r\prod_{i=1}^n{\Gamma(a+k(i-1))\Gamma(b)\over\Gamma(a+b+k-2n-i-1)}
\times\\
\times\prod_{i=1}^{n-r}(a+k(i-1))\prod_{i=n-r+1}^n(a+b+k(2n-i-1)){\goth
P}_r(a,b)\nonumber
\end{eqnarray}
where ${\goth P}_r(a,b)$ is a polynomial in $a$ and $b$.

On the other hand, we expand the hyperdeterminantal representation
of Aomoto's
 integral
by (\ref{devhank}) and obtain, after extracting common factors as
above,
\begin{eqnarray}
\Det({\bf A}_n^{a,b;k})&=&\Det({\bf S}_n(a,b;k))\sum_\lambda
{c_\lambda^{n,k}\over
Q(b)}\prod_{i=1}^n\left\{y(a)_{k(i-1)}^{\lambda_i-1}
(a+b)_{\lambda_i}^{k(n+i-2)-1}\right.\nonumber\\
&&\left.-(a)_{k(i-1)}^{\lambda_i}(a+b)_{\lambda_i+1}^{k(n+i-2)-1}\right\}
\end{eqnarray}
where $Q(b)$ is a polynomial of degree $kn(n-1)$ in $b$, $(a)_n^m$
denotes $(a+n)(a+n+1)\cdots(a+m)$ if $n\leq m$ and $(a)_n^m=1$
otherwise and the sum is over the partitions
$\lambda=(\lambda_1\leq \cdots\leq \lambda_n)$ verifying
$k(i-1)\leq \lambda_i\leq k(n+i-1)$ for each $i\in\{1,\cdots,n\}$.
After combining the coefficients of the $y^i$'s, one obtains
\begin{equation}\label{Ao2}
\Det({\bf A}^{a,b;k}_n(y))=\Det({\bf
S}_n(a,b;k))\sum_{r=0}^ny^r{{\goth Q}_r(a,b)\over Q(b)}
\end{equation}
where ${\goth Q}_r(a,b)$ is a polynomial.
 Hence, comparing the coefficients
 of  $y^i$ in the expressions (\ref{Ao1}) and (\ref{Ao2}), we find
\begin{eqnarray}\label{devaomoto}
\Det({\bf
A}^{a,b;k}_n(y))=\Det({\bf S}_n(a,b;k))\prod_{i=0}^{n-1}{1\over a+b+k(n+i-1)}\times\\\times
\sum_{i=0}^ny^i{P_i(a,b)\over
Q(b)}\prod_{j=1}^{n-i}(a+k(n-j))\prod_{j=n-i+1}^n(a+b+k(2n-j-1))\nonumber
\end{eqnarray}
where $P_i(a,b)$ is a polynomial. Remark that each $P_i(a,b)$
is of degree in $b$ at most $kn(n-1)$. Now, we apply (\ref{chvar})
and we equate the coefficients of  $y^i$ in
the left and right hand sides. After simplification, one obtains
for each $i$
\begin{eqnarray}\label{devaom2}
\sum_{j=i}^n\left(j\atop
i\right)P_j(b,a)F_j(b,a)Q(b)=
(-1)^{n-i}P_i(a,b)F_i(a,b)Q(a))
\end{eqnarray}
where
\begin{equation}\label{defF}F_i(a,b)=\displaystyle\prod_{m=1}^{n-i}(a+k(n-m))
\prod_{m=n-i+1}^{n}(a+b+k(2n-m-1)).\end{equation} Remarking that
$Q(b)$ does not divide $F_i(a,b)$, equality (\ref{devaom2})
implies
\begin{equation}
P_i(a,b)=R_i(a)Q(b).
\end{equation}
Hence, setting as in the previous section, $a=-k(n-1)$ in
\begin{eqnarray}
\sum_{j=i}^n\left(j\atop
i\right)R_j(b)F_j(b,a)=
(-1)^{n-i}R_i(a)F_i(a,b)
\end{eqnarray}
and remarking that
\begin{equation}F_i(-k(n-1),b)=0\end{equation} and
\begin{equation}F_i(b,-k(n-1))=\prod_{m=1}^n(b+k(n-m)),\end{equation}
  one obtains the recurrence relation
\begin{equation}
R_i(b)=-\sum_{j=i+1}^n\left(j\atop i\right)R_j(b).
\end{equation}
The starting point of the induction is $R_n(b)={1\over n!}$, which follows from  the fact that the coefficient of $y^n$
 in the expansion of the
hyperdeterminant is the value of  Selberg's integral.
Solving the recurrence relation,
 one finds that  each $R_i(b)$ is in fact independent of $b$,
and that
\begin{equation}
R_i(b)= {(-1)^{n-i}\over n!}\left(n\atop i\right)
\end{equation}
Hence,
\begin{eqnarray}
\Det({\bf A}_n^{a,b;k}(y))=
\Det({\bf S}_n(a,b;k)){1\over n!}\prod_{i=0}^{n-1}{1\over a+b+k(n+i-1)}\times\nonumber\\
\times\sum_{i=0}^n(-1)^{n-i}\left(n\atop
i\right)\prod_{m=1}^{n-i}(a+k(n-m))\prod_{m=n-i+1}^{n}(a+b+k(2n-m-1))y^i\\
=(-1)^n\Det({\bf S}_n(a,b;k)){\left({a\over k}\right)_n\over n!\left({a\over k}+{b\over k}+n-1\right)_n}\sum_{i=0}^n{(-n)_i({a\over k}+{b\over
k}+n-1)_i\over ({a\over k})_i}{y^i\over i!}\nonumber.\end{eqnarray}
We recognize a hypergeometric function of type $_2F_1$ which can
be evaluated as a monic Jacobi polynomial (see \cite{AAR} for
details)
\begin{eqnarray}
\Det({\bf A}_n^{a,b;k}(y))&=&(-1)^n\Det({\bf S}_n(a,b;k))
{\left({a\over k}\right)_n\over n!\left({a\over k}+{b\over
k}+n-1\right)_n}\times\nonumber\\ && \times _2F_1\left.\left(-n,
n+{a\over k}+{b\over k}-1\atop {a\over k}\right|y\right)\\
&=&{(-2)^{-n}\over n!}\Det({\bf S}_n(a,b;k))P_n^{{a\over k}-1,
{b\over k}-1}(1-2y)\nonumber
\end{eqnarray}
Aomoto's equality is therefore proved.

\section{Conclusion}

The examples discussed in this article show that  the
hyperdeterminantal calculus is a
pertinent tool to handle an interesting class of multiple integrals.
The hyperdeterminant is a particular invariant of hypermatrices. In
the case of  antisymmetric hypermatrices, another
invariant has similar properties: the hyperpfaffian
\cite{Bar,LT1}. In the classical case of matrices, de Bruijn has shown in
\cite{dB} that multiple integrals of some determinants could be
evaluated as Pfaffians. In \cite{LT1}, we have found some
generalizations to hyperpfaffians of these identities. In fact,
the generalized Heine theorem is a particular case of the
generalized de Bruijn integral, where the determinant can be
factored into a product of determinants of matrices whose
dimension is the number of integration variables. There
exist other  invariant polynomials of hypermatrices
\cite{Ca0,Ca1,Ca2,GKZ} which can be, in principle, computed using
Cayley's Omega process
or other methods of invariant theory (see \cite{LT3,BLT,BLTV} for examples
involving $2\times2\times2\times2$ and $3\times3\times3$
hypermatrices). A natural question is whether  there exist other
integral identities involving some of them.

Other generalizations of Selberg integral are
encountered in the  physical litterature. For example
the Dotsenko-Fateev (see \cite{DF}) and  Kaneko integrals (see
\cite{Kan}) give the partition functions of various systems.

Kaneko's integral reads
\begin{equation}
K_{n,\gamma}^{a,b}(y_1,\cdots,y_m)=\int_0^1|\Delta(x)|^{2\gamma}\prod_{i=1}^n
\left(x_i^{a-1}(1-x_i)^{b-1}\prod_{j=1}^n(x_i-y_j)dx_i\right)
\end{equation}
and is evaluated in terms of generalized orthogonal polynomials associated with Jack
polynomials (see \cite{Las1,Las2,Las3}). By
the generalized Heine theorem, when $\gamma$ is an integer,  Kaneko's integral
can be evaluated as a hyperdeterminant of moments (see \cite{LT2}).
Kaneko's proof of his identity
is related to the Calogero-Sutherland model. Up to now, we have been unable
unable to find a hyperdeterminantal interpretation of this. In the
same way, Anderson's (see \cite{An}) and Aomoto's (see \cite{Ao1})
proofs of Selberg's integral
seem to contain  information unrelated to the hyperdeterminantal
representations.

\end{document}